\pdfoutput=1
\RequirePackage{ifpdf}
\ifpdf 
\documentclass[pdftex]{sigma}
\else
\documentclass{sigma}
\fi

\numberwithin{equation}{section}

\newtheorem{Theorem}{Theorem}[section]
\newtheorem*{Theorem*}{Theorem}
\newtheorem{Corollary}[Theorem]{Corollary}
\newtheorem{Lemma}[Theorem]{Lemma}
\newtheorem{Proposition}[Theorem]{Proposition}

\theoremstyle{definition}

\newtheorem*{Example*}{Example}

\newtheorem*{Remark*}{Remark}
\newtheorem*{Remarks*}{Remarks}

\newcommand{\R}{{\mathbb R}}
\newcommand{\IR}{{\mathbb R}}
\newcommand{\C}{{\mathbb C}}
\newcommand{\IC}{{\mathbb C}}
\newcommand{\Z}{\mathbb{Z}}

\newcommand{\IN}{{\mathbb N}}
\newcommand{\GL}{\mathrm{GL}}
\newcommand{\Aut}{\operatorname{Aut}}
\newcommand{\Mat}{\mathrm{Mat}}
\newcommand{\Nilp}{\mathrm{Nilp}}

\begin{document}
\allowdisplaybreaks

\newcommand{\arXivNumber}{2509.16178}

\renewcommand{\thefootnote}{}

\renewcommand{\PaperNumber}{016}

\FirstPageHeading

\ShortArticleName{Asymptotics for the Enumeration of Commuting Matrices over Finite Fields}

\ArticleName{Asymptotics for the Enumeration\\ of Commuting Matrices over Finite Fields\footnote{This paper is a~contribution to the Special Issue on Recent Advances in Vertex Operator Algebras in honor of James Lepowsky. The~full collection is available at \href{https://sigma-journal.com/Lepowsky.html}{https://sigma-journal.com/Lepowsky.html}}}

\Author{Kathrin BRINGMANN~$^{\rm a}$, Shane CHERN~$^{\rm b}$, Johann FRANKE~$^{\rm a}$ and Bernhard HEIM~$^{\rm a}$}

\AuthorNameForHeading{K.~Bringmann, S.~Chern, J.~Franke and B.~Heim}

\Address{$^{\rm a)}$~University of Cologne, Department of Mathematics and Computer Science,\\
\hphantom{$^{\rm a)}$}~Weyertal 86-90, 50931 Cologne, Germany}
\EmailD{\mail{kbringma@uni-koeln.de}, \mail{jfrank12@uni-koeln.de}, \mail{bheim@uni-koeln.de}}

\Address{$^{\rm b)}$~Fakult\"at f\"ur Mathematik, Universit\"at Wien, Oskar-Morgenstern-Platz 1, Wien 1090, Austria}
\EmailD{\mail{xiaohangc92@univie.ac.at}}

\ArticleDates{Received September 22, 2025, in final form January 19, 2026; Published online February 19, 2026}

\Abstract{We give asymptotic expressions for the number of commuting matrices over finite fields. For this, we use product expansions for the corresponding generating functions.}


\Classification{11N45; 05A16}

\begin{flushright}
{\it In honor of James Lepowsky's 80-th birthday}
\end{flushright}

\renewcommand{\thefootnote}{\arabic{footnote}}
\setcounter{footnote}{0}

\section{Introduction and statement of results}

Let $\mathbb{F}_{p^r}$ be the finite field with $p^r$ elements, where $p$ is a prime and $r\in \mathbb{N}$. The study of the enumeration of commuting square matrix pairs over $\mathbb{F}_{p^r}$ was first launched by Feit and Fine~\cite{FF}. Let $Q_{p^r}(n)$ denote the number of ordered pairs of (not necessarily distinct) commuting $n \times n$ matrices with entries in $\mathbb{F}_{p^r}$. Feit and Fine proved that $Q_{p^r}(n)$ (re-normalized) admits the following product-like generating function.

\begin{Theorem}[Feit--Fine, {\cite[p.~91]{FF}}]
	We have
	\begin{align} \label{eq:Power-series}
		\sum_{n \geq 0} \frac{Q_{p^r}(n)}{p^{rn^2}f_{p^r}(n) } w^n = \prod_{\substack{\ell \geq 1 \\ j \geq 0}} \frac{1}{1-p^{r(1-j)}w^\ell},
	\end{align}
	where
	\begin{align*}
		f_{p^r}(n) := \prod_{j=1}^n \left( 1 - p^{-rj}\right).
	\end{align*}
\end{Theorem}

\begin{Remarks*}\quad
	\begin{enumerate}
		\item[$(1)$] The product in \eqref{eq:Power-series} converges for $w\in \mathbb{C}$ with $|w| < 1$ not being a pole and continues to a~meromorphic function for $|w| < 1$. This is due to the absolute and uniform convergence on compact subsets of the double series
		\[\sum_{\substack{\ell \geq 1 \\ j \geq 0}} p^{r(1-j)}w^\ell\]
		away from a pole.
		\item[$(2)$] Recall that (see \cite[equation~(2.1)]{FGS17})
		\begin{equation} \label{eq:GL-identity}
			\left| \mathrm{GL}_n(\mathbb{F}_{p^r}) \right|
			= p^{r n^2} f_{p^r}(n).
		\end{equation}
		Therefore, \eqref{eq:Power-series} can alternatively be written as
		\begin{align*}
			\sum_{n \geq 0} \frac{Q_{p^r}(n)}{	\left| \mathrm{GL}_n(\mathbb{F}_{p^r}) \right|} w^n = \prod_{\substack{\ell \geq 1 \\ j \geq 0}} \frac{1}{1-p^{r(1-j)}w^\ell}.
		\end{align*}
		\item[$(3)$] By \cite[equation~(1)]{FF} and \eqref{eq:GL-identity}, we have
		\begin{equation*}
			Q_{p^r}(n) = \vert \mathrm{GL}_n(\mathbb{F}_{p^r})\vert
			\sum_{\lambda \vdash n}
			\frac{p^{r \sum_{k=1}^{n} b_k}}{\prod_{k=1}^{n}{f_{p^r}(b_k)}},
		\end{equation*}
		where the sum runs over all partitions of $n$ and $\lambda$ has the frequencies $b_k$ given by $n=\sum_{k=1}^n b_k k$.
	\end{enumerate}
\end{Remarks*}

A natural question is to find an approximation for $Q_{p^r}(n)$, as $n \to \infty$. It was shown by Motzkin--Taussky~\cite{Mot} (see also \cite[p.~71, Theorem~2]{Gur1992}) that the space of $n\times n$ commuting matrices over an algebraically closed field has dimension $n^2+n$. Hence, heuristically, $Q_{p^r}(n)$ should be proportional to \smash{$p^{r(n^2+n)}$} for $n$ sufficiently large. This observation was rigorously proved by Fulman and Guralnick~\cite[Theorem~2.6]{FG18} as a byproduct of their new proof of the result of Feit and Fine. More precisely, Fulman and Guralnick established that
\begin{equation*} 
	\lim_{n \to \infty} \frac{Q_{p^r}(n)}{p^{r(n^2+n)}} = \prod_{j\ge1} \left( 1 - p^{-rj}\right)^{-j}.
\end{equation*}

In this paper, we prove more precise asymptotics for $Q_{p^r}(n)$.
We start by writing for $w \in \IC$ with $|w| < 1$,
\begin{align} \label{eq:prod-zerlegung}
	\prod_{\substack{\ell \geq 1 \\ j \geq 0}} \frac{1}{1-p^{r(1-j)}w^\ell} = P_{0, p^r}(w) F_{p^r}(w),
\end{align}
where for $x>0$,
\begin{align} \label{eq:P-n-def}
	P_{m, x}(w) := \prod_{\substack{\ell \geq 1 \\\ell \not= m} }\frac{1}{1-xw^\ell}, \qquad F_{p^r}(w) := \prod_{\ell,j \geq 1} \frac{1}{1-p^{r(1-j)}w^\ell}.
\end{align}

Theorem~\ref{T:QAsymp} gives the asymptotic behavior of $Q_{p^r}(n)$ as $n$ gets large. In it, and indeed throughout this note, we use the notation \smash{$\zeta_m := {\rm e}^{\frac{2\pi {\rm i}}{m}}$}.

\begin{Theorem}\label{T:QAsymp}
	We have, as $n\to\infty$,
	\begin{align} \label{eq:Q-expansion}
		Q_{p^r}(n) \sim p^{rn^2} f_{p^r}(n) \sum_{m \geq 1} C_{m,p^r}(n) p^{\frac{rn}{m}},
	\end{align}
	where
	\begin{align*}
		C_{m,p^r}(n) := \frac{1}{m} \sum_{j=0}^{m-1} P_{m,p^r}\left(\zeta_m^{-j} p^{-\frac{r}{m}}\right)F_{p^r}\left(\zeta_m^{-j} p^{-\frac{r}{m}}\right) \zeta_m^{nj}.
	\end{align*}
	More precisely, we have, as $n\to \infty$,
	\begin{align*}
		{Q_{p^r}(n) = p^{r(n^2+n)} \prod_{j \geq 1} \left( 1 - p^{-rj}\right)^{-j}} + O\left(p^{rn^2 + \frac{rn}{2}} \right).
	\end{align*}	
\end{Theorem}

\begin{Remark*}
	The series in \eqref{eq:Q-expansion} does not converge, but should be viewed as an asymptotic expansion. Also note that $C_{m,p^r}(n)$ only depends on $n\pmod{m}$. In particular, $C_{1,p^r}(n)$ does not depend on $n$. Furthermore, we have
	\begin{align*}
		|C_{m,p^r}(n)| \leq \max_{0 \leq j \le m-1} \left| P_{m,p^r}\left(\zeta_m^{-j} p^{-\frac{r}{m}}\right)F_{p^r}\left(\zeta_m^{-j} p^{-\frac{r}{m}}\right)\right| = \left|P_{m,p^r}\left(p^{-\frac{r}{m}}\right)F_{p^r}\left(p^{-\frac{r}{m}}\right) \right|.
	\end{align*}
	As a result, when truncated at $m=N$, the error in \eqref{eq:Q-expansion} is $O\left(p^{\frac{rn}{N+1}}\right)$.
\end{Remark*}

\begin{Example*}
	Let $p=2$ and $r=1$. By \eqref{eq:Power-series}, we have
	\begin{align} \label{eq:1-iteration}
		(1 - 2w) \sum_{n \geq 0} \frac{Q_{2}(n)}{2^{n^2}f_{2}(n) } w^n = \prod_{\substack{\ell \geq 1 \\ j \geq 0 \\ (\ell, j) \not= (1,0)}} \frac{1}{1-2^{1-j}w^\ell}.
	\end{align}
	While \eqref{eq:Power-series} has radius of convergence \smash{$\tfrac12$}, the right-hand side of \eqref{eq:1-iteration} has radius of convergence~\smash{$\frac{1}{\sqrt{2}}$}. Hence, we obtain
	\begin{align*}
		C_{1,2}(1) = \lim_{n \to \infty} \frac{Q_{2}(n)}{2^{n^2+n}f_2(n)} = \prod_{j \geq 1} \frac{1}{\left( 1 - 2^{-j}\right)^{j+1}} = 34.738723457 \ldots.
	\end{align*}
	Further values are given by
	\begin{alignat*}{3}
		&C_{2,2}(1) = -11716.7651425569 \ldots, \qquad&& C_{2,2}(2) = -11716.3960075313\ldots,& \\
		&C_{3,2}(1) = 7970793.64416118\ldots, \qquad&& C_{3,2}(2) = 7970793.59033743 \ldots,& \\
		&C_{3,2}(3) = 7970793.67801128 \ldots .&
	\end{alignat*}
	Looking at these first values, one notes that they are very close to each other for $m$ fixed. A~further investigation of their behavior would be interesting.
\end{Example*}


The paper is organized as follows. First we recall some basic definitions. Then in Section~\ref{sec:QAsymp}, we give a proof of Theorem~\ref{T:QAsymp}. In Section~\ref{sec:CL-series}, we introduce the concept of Cohen--Lenstra series, including the series in \eqref{eq:Power-series} as an instance. We generalize the ideas in Section~\ref{sec:QAsymp} to another Cohen--Lenstra series in Theorem~\ref{thm:formula-2}. In the final section, we propose problems for further research.

\section{Preliminaries}

In the following, we recall some required results from complex analysis. We start with a well-known and simple observation that is still very useful for many aspects in the theory of $q$-series, recurrences, and combinatorics. For $\varrho\in \R^+$, let $B_\varrho(a) := \{z\in\mathbb{C} : |z-a|<\varrho\}$.

\begin{Proposition} 
	Let $a \in \IC$, $\varrho \in \IR^+$, and $f$ be a holomorphic function defined in a neighborhood of~$a$. Then the following are equivalent:
	\begin{enumerate}\itemsep=0pt
		\item[$(1)$] The function $f$ extends to an analytic function on $B_\varrho(a)$, but not on $B_{\varrho+\varepsilon}(a)$ for any $\varepsilon > 0$.
		\item[$(2)$] The power series \smash{$\sum_{n\geq 0} \frac{f^{(n)}(a)}{n!} (z-a)^n$} has radius of convergence $\varrho$.
	\end{enumerate}
\end{Proposition}

For $w \in \IC$ and a sequence $a_m$ of non-zero complex numbers, such that $\sum_{m \geq 1}|a_m|^{-1} < \infty$, we let
\begin{align}\label{F}
	F(w) := \prod_{m \geq 1} \left( 1 - \frac{w}{a_m}\right).
\end{align}

\begin{Proposition}[Wang, {\cite[Theorem~1]{Wang}}] \label{prop:product}
	Let $a_m$ be a sequence of non-zero complex numbers such that $a_m \not= a_n$ for $m \not= n$ and such that \smash{$\sum_{m \geq 1} |a_m|^{-1}$} converges. Then the function \smash{$\frac{1}{F}$} is holomorphic in a neighborhood of $w=0$ and we have 
	\begin{align*} 
		\frac{1}{F(w)} = -\sum_{k \geq 0} \sum_{m \geq 1} \frac{1}{F'(a_m)a_{m}^{k+1}} w^k.
	\end{align*}
\end{Proposition}

\section{Proof of Theorem~\ref{T:QAsymp}}\label{sec:QAsymp}

To prove Theorem~\ref{T:QAsymp}, we require the following lemma.

\begin{Lemma} \label{lem:Removable}
	Let $x > 1$, $n \in \mathbb{N}$, \smash{$0 < \varepsilon < x^{-\frac1n}$}, and $f$ be a holomorphic function in some region \smash{$\mathcal{R}_{x,n,\varepsilon}:=\{ w \in \C  :  x^{-\frac1n} - \varepsilon < |w| < x^{-\frac1n} + \varepsilon\}$}. Define
	\begin{align*}
		g(w) := \frac{f(w)}{1-x w^n}.
	\end{align*}
	Then
	\begin{align*}
		g(w) - \frac{1}{n}\sum_{j=0}^{n-1} \frac{f\left(\zeta_n^{-j} x^{-\frac1n}\right)}{1 - \zeta_n^j x^\frac1n w}
	\end{align*}
	continues holomorphically to $\mathcal{R}_{x,n,\varepsilon}$.
\end{Lemma}

\begin{proof}
	The only poles of $g$ in $\mathcal{R}_{x,n,\varepsilon}$ lie in $\{ \zeta_n^{-j} x^{-\frac{1}{n}}  :  0 \leq j \le n-1\}$. These are simple. We compute, using l'Hospital,
	\begin{align*}
		\lim_{w \to \zeta_n^{-j} x^{-\frac1n}} \frac{\left(1 - \zeta_n^j x^\frac1n w\right) f(w)}{1-x w^n} = \frac{1}{n} f\left(\zeta_n^{-j} x^{-\frac1n}\right).
	\end{align*}
	The claim follows.
\end{proof}

We have the following corollary to Lemma~\ref{lem:Removable}.

\begin{Corollary} \label{cor:poles}
	Let $f$ be a function holomorphic inside the unit disk and $N \in \mathbb{N}$. Then
	\begin{align*}
		P_{0,p^r}(w) f(w) - \sum_{m=1}^N \frac{1}{m} \sum_{j=0}^{m-1} \frac{P_{m,p^r}\left(\zeta_m^{-j} p^{-\frac{r}{m}}\right)f\left(\zeta_m^{-j} p^{-\frac{r}{m}}\right)}{1 - \zeta_m^j p^\frac{r}{m} w}
	\end{align*}
	continues holomorphically to the region $\{w\in\C  :  |w| < p^{-\frac{r}{N+1}}\}$.
\end{Corollary}

\begin{proof}
	By \eqref{eq:P-n-def}, we have, for $m\in\IN$,
	\begin{align*}
		P_{0,p^r}(w)f(w) = \frac{P_{m,p^r}(w)f(w)}{1 - p^r w^m}.
	\end{align*}
	Note that the numerator is holomorphic in the region
	\begin{align*}
		\left\{ w \in \C  :  p^{-\frac{r}{m-1}} < |w| < p^{-\frac{r}{m+1}} \right\},
	\end{align*}
	where \smash{$p^{-\frac{r}{0}} := 0$}. By applying Lemma~\ref{lem:Removable} iteratively, the claim follows.
\end{proof}

We are now ready to prove Theorem~\ref{T:QAsymp}.

\begin{proof}[Proof of Theorem~\ref{T:QAsymp}]
	Let $N \in \IN$ be arbitrary. We abbreviate
	\begin{align*}
		c_{m,j} := P_{m,p^r}\left(\zeta_m^{-j} p^{-\frac{r}{m}}\right)F_{p^r}\left(\zeta_m^{-j} p^{-\frac{r}{m}}\right).
	\end{align*}
	By \eqref{eq:Power-series} and \eqref{eq:prod-zerlegung}, and Corollary~\ref{cor:poles}, we conclude that
	\begin{align*}
		\sum_{n \geq 0} \frac{Q_{p^r}(n)}{p^{rn^2}f_{p^r}(n) } w^n - \sum_{m=1}^N \frac{1}{m} \sum_{j=0}^{m-1} \frac{c_{m,j}}{1 - \zeta_m^j p^\frac{r}{m} w}
	\end{align*}
	continues holomorphically to $\{w\in\C  :  |w| < p^{-\frac{r}{N+1}}\}$. We can then write
	\begin{align*}
		\sum_{m=1}^N \frac{1}{m} \sum_{j=0}^{m-1} \frac{c_{m,j}}{1 - \zeta_m^j p^\frac{r}{m} w} = \sum_{n \geq 0} \sum_{m=1}^N \frac{p^\frac{rn}{m} w^n}{m} \sum_{j=0}^{m-1} c_{m,j} \zeta_m^{jn}.
	\end{align*}
	Hence we have, for $\varepsilon > 0$,
	\begin{align*}
		\frac{Q_{p^r}(n)}{p^{rn^2}f_{p^r}(n) } - \sum_{m=1}^N \frac{1}{m} \sum_{j=0}^{m-1} c_{m,j} \zeta_m^{jn} p^\frac{rn}{m} = O_\varepsilon \left( p^{\left(\frac{r}{N+1} + \varepsilon \right)n}\right), \qquad n \to \infty.
	\end{align*}
	Equivalently, we have, as $n\to\infty$,
	\begin{equation*}
		Q_{p^r}(n)= p^{rn^2}f_{p^r}(n) \sum_{m=1}^N \frac{1}{m} \sum_{j=0}^{m-1} c_{m,j} \zeta_m^{jn} p^{\frac{rn}{m}} + O_\varepsilon \left( p^{rn^2 + \left(\frac{r}{N+1} + \varepsilon \right)n}\right).
	\end{equation*}
	The result follows after a straightforward calculation.
\end{proof}

Finally, we give an upper bound for the coefficients in the asymptotic expansion in Theorem~\ref{T:QAsymp}.

\begin{Proposition} 
	We have, for $\varepsilon > 0$, uniformly in $m$, $p^r$, and $n$,
	\begin{align*}
		C_{m,p^r}(n) \ll_\varepsilon \exp\left( \frac{\zeta(3) + \varepsilon }{\left(1 - p^{-\frac{r}{m}}\right)^2}\right),
	\end{align*}
	where $\zeta$ denotes the Riemann zeta function.
\end{Proposition}

\begin{proof}
	By \eqref{eq:P-n-def}, we have
	\begin{align*}
		\left|P_{m,p^r}\left({{p}^{-\frac{r}{m}}}\right)F_{p^r}\left({{p}^{-\frac{r}{m}}}\right) \right| = \prod_{\substack{\ell \geq 1 \\\ell \not= m}} \left| \frac{1}{1-p^{r - \frac{r\ell}{m}}} \right| \prod_{\ell,j \geq 1} \left| \frac{1}{1-p^{r(1-j)}p^{-\frac{r\ell}{m}}} \right|.
	\end{align*}
	We estimate both products separately. First, we have
	\begin{align*}
		\prod_{\substack{\ell \geq 1 \\\ell \not= m}} \frac{1}{1-p^{r - \frac{r\ell}{m}}} = (-1)^{m+1} p^{-\frac{rm(m+1)}{2m}} P_{0,1}\left(p^{-\frac{r}{m}}\right) \prod_{\ell=1}^{m-1} \frac{1}{1 - p^{-\frac{r\ell}{m}}} \ll p^{-\frac{r(m-1)}{2}} P_{0,1}^2\left(p^{-\frac{r}{m}}\right).
	\end{align*}
	For the second product, we find
	\begin{equation*}
		\prod_{\ell,j \geq 1} \frac{1}{1-p^{r(1-j)}p^{-\frac{r\ell}{m}}} = \prod_{\substack{j \geq 0\\ \ell \geq 1}} \frac{1}{1 - p^{-\frac{r}{m}(mj+\ell)}} \ll \prod_{\substack{j \geq 0\\ \ell \geq 1}} \frac{1}{1 - p^{-\frac{r}{m}(j+\ell)}} = \prod_{j \geq 1} \left(1 - p^{-\frac{rj}{m}}\right)^{-j}.
	\end{equation*}
	Applying \cite[Lemma~3.2]{BBBF23}, we have, as $w \to 1^-$,
	\begin{align*}
		P_{0,1}^2\left( w \right) \ll_\varepsilon {\rm e}^{\frac{\zeta(2)+ \varepsilon}{|w-1|}}, \qquad \prod_{j \geq 1} \frac{1}{(1 - w^j)^j} \ll_\varepsilon {\rm e}^{\frac{\zeta(3)+ \varepsilon}{|w-1|^2}}
	\end{align*}
	for all $\varepsilon > 0$. As a result, for all $\varepsilon > 0$, we find the claimed bound.
\end{proof}

\section{Cohen--Lenstra series and the nilpotent classes}\label{sec:CL-series}

Since its discovery, Feit and Fine's generating function \eqref{eq:Power-series} has been reproved by other means such as the motivic Donaldson--Thomas theory \cite{BM2015}. More recently, Huang~\cite{Hua2023} placed such counting problems into
the framework of Cohen--Lenstra series \cite{CL1984}. Simply speaking, Huang considered \emph{$($generalized$)$ Cohen--Lenstra series}, which are defined as
\begin{align*}
	Z_{R}(w) = Z_{R/\mathbb{F}_{p^r}}(w) := \sum_{M} \frac{1}{|{\Aut(M)}|} w^{\dim_{\mathbb{F}_{p^r}} (M)},
\end{align*}
where $M$ runs over all isomorphism classes of finite $R$-modules and $R$ is a commutative ring containing the finite field $\mathbb{F}_{p^r}$ such that the above series is well-defined.\footnote{The original Cohen--Lenstra series in \cite{CL1984} are defined for $R$ a Dedekind domain over $\mathbb{F}_{p^r}$.} A more general setting was later provided by Huang in \cite{Hua2024}. Note that, by \eqref{eq:GL-identity}, the relation \eqref{eq:Power-series} corresponds to the Cohen--Lenstra series (see also Huang's list in \cite[p.~27, item 3]{Hua2023}):
\begin{align*}
	Z_{\mathbb{F}_{p^r}[u,v]}(w) = \sum_{n\ge 0} \frac{\left|\left\{(A,B)\in \Mat_n(\mathbb{F}_{p^r})^2  :  AB=BA\right\}\right|}{|\GL_n(\mathbb{F}_{p^r})|} w^n = \prod_{\substack{\ell \geq 1 \\ j \geq 0}} \frac{1}{1-{p}^{r(1-j)}w^\ell},
\end{align*}
where, as usual, $\Mat_n(\mathbb{F}_{p^r})$ denotes the set of $n\times n$ matrices over $\mathbb{F}_{p^r}$. In this setting, we may generalize $R$ to count commuting varieties (i.e., subvarieties of $XY - YX = 0$) over a finite field; see \cite[p.~40, Proposition~4.3]{Hua2023}. However, these Cohen--Lenstra series in general do not admit a product form and mostly behave multisum-like; see \cite[equations~(1.14) and (1.15)]{Hua2023}. Therefore, our method of proof used for Theorem~\ref{T:QAsymp} does not apply. To work out new methods to deal with these series would be worthwhile.
Failing this, we restrict to the easier problem of enumerating commuting square matrix pairs in nilpotent classes. Let
\begin{align*}
	\Nilp_n(\mathbb{F}_{p^r}):=\left\{A\in\Mat_n(\mathbb{F}_{p^r}) :  \text{$A$ is nilpotent}\right\}.
\end{align*}
In the simplest cases, due to Fine and Herstein~\cite{FH1958} and Fulman and Guralnick~\cite{FG18}, we still have product-like generating functions similar to \eqref{eq:Power-series}. In terms of the Cohen--Lenstra series, the two cases may be written as follows:

\begin{itemize}\itemsep=0pt
	\item $R=\mathbb{F}_{p^r}[[u]]$ -- Fine--Herstein~\cite[p.~499, Theorem~1]{FH1958} (see also the list of Huang \cite[p.~27, item~2]{Hua2023}):
	\begin{align} \label{eq:Nil-identity}
		Z_{\mathbb{F}_{p^r}[[u]]}(w) = \sum_{n\ge 0} \frac{|\Nilp_n(\mathbb{F}_{p^r})|}{|\GL_n(\mathbb{F}_{p^r})|} w^n = \prod_{j\ge 1} \frac{1}{1-p^{-rj}w},
	\end{align}
	
	\item $R=\mathbb{F}_{p^r}[[u,v]]$ -- Fulman--Guralnick~\cite[p.~301, Theorem~2.9]{FG18} (see also the list of Huang~\cite[p.~27, item~4]{Hua2023}):
	\begin{align}
		Z_{\mathbb{F}_{p^r}[[u,v]]}(w) & = \sum_{n\ge 0} \frac{\left|\left\{(A,B)\in \Nilp_n(\mathbb{F}_{p^r})^2 :  AB=BA\right\}\right|}{|\GL_n(\mathbb{F}_{p^r})|} w^n \nonumber\\
& = \prod_{\substack{\ell \geq 1 \\ j \geq 0}} \frac{1}{1- p^{-r(1+j)}w^\ell}.\label{eq:Nil-identity-2}
	\end{align}
\end{itemize}

The structure of \eqref{eq:Nil-identity} is particularly simple, which allows us to use Proposition~\ref{prop:product} to derive not only an asymptotic expression but even a closed-form series representation for its coefficients. For $w \in \IC$, we define
\begin{align*}
	Z_{m, \mathbb{F}_{p^r}[[u]]}(w) := \prod_{\substack{ j\ge 1 \\ j \not= m}} \frac{1}{1-p^{-rj}w}.
\end{align*}

\begin{Theorem} \label{thm:formula-2}
	Let $p$ be a prime and $r \in \mathbb{N}$. Then we have, for $n \in \mathbb{N}_0$,
	\begin{align*}
		\frac{|\Nilp_n(\mathbb{F}_{p^r})|}{|\GL_n(\mathbb{F}_{p^r})|} = \sum_{m \geq 1} \frac{Z_{m, \mathbb{F}_{p^r}[[u]]}(p^{rm})}{ p^{rnm}}.
	\end{align*}
\end{Theorem}

\begin{proof}
	We apply Proposition~\ref{prop:product} with $a_j := p^{rj}$. As $p^r > 1$, the series $\sum_{m \geq 1} p^{-rm}$ converges and $a_m \not= a_n$ for all $m \not= n$. By the product rule (see \eqref{F} for the definition of $F$),
	\begin{align*} 
		F'(p^{rm}) = - \frac{1}{p^{rm} Z_{m, \mathbb{F}_{p^r}[[u]]}(p^{rm})}.
	\end{align*}	
	Using this, 
	we obtain the theorem.
\end{proof}

\begin{Example*}
	The series in Theorem~\ref{thm:formula-2} is actually convergent. We choose $p=2$ and $r=1$ in this example. And let
	\begin{align*}
		Z_{m, \mathbb{F}_{2}[[u]], N}(w) := {\prod_{\substack{N \ge j\ge 1 \\ j \not= m}} \frac{1}{1- 2^{-j}w}.}
	\end{align*}
	We have
	\begin{align*}
		Z_{\mathbb{F}_{2}[[u]]}(w) = \prod_{j\ge 1} \frac{1}{1-2^{-j}w} = 1 + w + \frac{2}{3}w^2 + \frac{4}{11}w^3 + \frac{64}{315} w^4 + \cdots.
	\end{align*}
	On the other hand, for $N=100$, one finds numerically
	\begin{align*}
		&\sum_{1\le m\le 10} Z_{m, \mathbb{F}_{2}[[u]],100}(2^m)= 0.999999999999999667 \ldots, \\
		&\sum_{1\le m\le 10} \frac{Z_{m, \mathbb{F}_{2}[[u]],100}(2^m)}{2^{m}}= 0.9999999999999999998 \ldots, \\
		&\sum_{1\le m\le 10} \frac{Z_{m, \mathbb{F}_{2}[[u]],100}(2^m)}{2^{2m}}= 0.6666666666666666666 \ldots, \\
		&\sum_{1\le m\le 10} \frac{Z_{m, \mathbb{F}_{2}[[u]],100}(2^m)}{2^{3m}}= 0.3809523809523809523 \ldots, \\
		&\sum_{1\le m\le 10} \frac{Z_{m, \mathbb{F}_{2}[[u]],100}(2^m)}{2^{4m}}= 0.2031746031746031746 \ldots.
	\end{align*}
	These sums represent the cases $0\leq n \leq 4$. Note that also the degenerate case $n=0$ is included.
\end{Example*}

\section{Questions for future research}

\noindent We propose the following research questions.
\begin{enumerate}\itemsep=0pt
	\item[$(1)$] Is there any recursive formula for the $Q_{p^r}(n)$ (or quotients of it by normalization factors) that is analogous to the divisor recurrence formula or the Euler pentagonal theorem \cite[equation~(1.3.1)]{And1998}
	\begin{align*}
		\prod_{n \geq 1} (1 - q^n) = \sum_{n \in \Z} (-1)^n q^{\frac{n(3n-1)}{2}}?
	\end{align*}
	In the second scenario, there is no exact identity, but rather an error term smaller than all exponential functions in the main series. Is there a combinatorial interpretation for these recurrences and the error term?
	
	\item[$(2)$] Can we improve the bounds for $C_{m,p^r}(n)$? What can we say about these coefficients for fixed $m$?
	
	\item[$(3)$] How could the methods in this paper be extended to deal with multivariable products such as \eqref{eq:Nil-identity-2}?
\end{enumerate}

\subsection*{Acknowledgements}

The first author has received funding from the European Research Council (ERC) under the European Union’s Horizon 2020 research and innovation programme (grant agreement No. 101001179). The second author has been supported by the Austrian Science Fund (grant agreement No. 10.55776/F1002). The authors are grateful to Yifeng Huang for helpful conversations regarding the background of Cohen--Lenstra series. Moreover, we thank the referees for their helpful comments.

\pdfbookmark[1]{References}{ref}
\LastPageEnding


\begin{thebibliography}{99}
\footnotesize\itemsep=0pt

\bibitem{And1998}
Andrews G.E., The theory of partitions, \textit{Cambridge Math. Lib.},
 \href{https://doi.org/10.1017/CBO9780511608650}{Cambridge University Press},
 Cambridge, 1998.

\bibitem{BBBF23}
Bridges W., Brindle B., Bringmann K., Franke J., Asymptotic expansions for
 partitions generated by infinite products,
 \href{https://doi.org/10.1007/s00208-024-02807-x}{\textit{Math. Ann.}}
 \textbf{390} (2024), 2593--2632,
 \href{http://arxiv.org/abs/2303.11864}{arXiv:2303.11864}.

\bibitem{BM2015}
Bryan J., Morrison A., Motivic classes of commuting varieties via power
 structures,
 \href{https://doi.org/10.1090/S1056-3911-2014-00657-3}{\textit{J.~Algebraic
 Geom.}} \textbf{24} (2015), 183--199,
 \href{http://arxiv.org/abs/1206.5864}{arXiv:1206.5864}.

\bibitem{CL1984}
Cohen H., Lenstra Jr. H.W., Heuristics on class groups of number fields, in
 Number Theory, {N}oordwijkerhout 1983 ({N}oordwijkerhout, 1983),
 \textit{Lecture Notes in Math.}, Vol.~1068,
 \href{https://doi.org/10.1007/BFb0099440}{Springer}, Berlin, 1984, 33--62.

\bibitem{FF}
Feit W., Fine N.J., Pairs of commuting matrices over a finite field,
 \href{https://doi.org/10.1215/S0012-7094-60-02709-5}{\textit{Duke Math.~J.}}
 \textbf{27} (1960), 91--94.

\bibitem{FH1958}
Fine N.J., Herstein I.N., The probability that a matrix be nilpotent,
 \href{https://doi.org/10.1215/ijm/1255454112}{\textit{Illinois~J. Math.}}
 \textbf{2} (1958), 499--504.

\bibitem{FG18}
Fulman J., Guralnick R., Enumeration of commuting pairs in {L}ie algebras over
 finite fields, \href{https://doi.org/10.1007/s00026-018-0390-4}{\textit{Ann.
 Comb.}} \textbf{22} (2018), 295--316,
 \href{http://arxiv.org/abs/1611.05499}{arXiv:1611.05499}.

\bibitem{FGS17}
Fulman J., Guralnick R., Stanton D., Asymptotics of the number of involutions
 in finite classical groups,
 \href{https://doi.org/10.1515/jgth-2017-0011}{\textit{J.~Group Theory}}
 \textbf{20} (2017), 871--902,
 \href{http://arxiv.org/abs/1602.03611}{arXiv:1602.03611}.

\bibitem{Gur1992}
Guralnick R.M., A~note on commuting pairs of matrices,
 \href{https://doi.org/10.1080/03081089208818123}{\textit{Linear Multilinear
 Algebra}} \textbf{31} (1992), 71--75.

\bibitem{Hua2024}
Huang Y., Commuting matrices via commuting endomorphisms,
 \href{http://arxiv.org/abs/2404.19483}{arXiv:2404.19483}.

\bibitem{Hua2023}
Huang Y., Mutually annihilating matrices, and a {C}ohen--{L}enstra series for
 the nodal singularity,
 \href{https://doi.org/10.1016/j.jalgebra.2022.11.021}{\textit{J.~Algebra}}
 \textbf{619} (2023), 26--50,
 \href{http://arxiv.org/abs/2110.15566}{arXiv:2110.15566}.

\bibitem{Mot}
Motzkin T.S., Taussky O., Pairs of matrices with property~{L}.~{II},
 \href{https://doi.org/10.2307/1992996}{\textit{Trans. Amer. Math. Soc.}}
 \textbf{80} (1955), 387--401.

\bibitem{Wang}
Wang X., The summation of infinite partial fraction decomposition~{I}: some
 formulae related to the {H}urwitz zeta function,
 \href{http://arxiv.org/abs/2102.04115}{arXiv:2102.04115}.

\end{thebibliography}
\end{document}